\documentclass[a4paper,12pt]{article}

\usepackage[top=4cm,bottom=3.4cm,left=2.4cm,right=3cm]{geometry}

\usepackage{textgreek}

\usepackage{arydshln}
\usepackage{amssymb,amsfonts,graphicx}

\newcommand{\Z}{\mathbb{Z}}

\makeatletter

\renewcommand{\@seccntformat}[1]{{\csname the#1\endcsname}{\normalsize .}\hspace{.5em}}
\makeatother

\usepackage{ifpdf}
\usepackage{indentfirst}

\def \[{\begin{equation}}
\def \]{\end{equation}}

\begin{document}

\setlength{\baselineskip}{15pt}
\begin{center}{\large \bf Table of large graphs with given degree and diameter }
\vspace{4mm}

{ Francesc Comellas \\

Departament de  Matem\`{a}tiques, Universitat Polit\`{e}cnica de Catalunya,\\
 Barcelona, Catalonia, Spain  \quad  {\tt\footnotesize  francesc.comellas@upc.edu}}

\end{center}

\noindent {\bf Abstract}:\,   We update the table of large undirected graphs with given degree and diameter with results obtained since the publication of the survey 
  by M. Miller and  J. \v{S}ir\'{a}\v{n} in the {\em Electronic Journal of Combinatorics} (Dynamic Survey DS14, 2nd. edition. May 2013).

\vspace{3mm} {\footnotesize  \noindent{\it Keywords: Large Graphs, Degree-Diameter Problem,  Cayley Graphs} }

\section{Introduction}

The construction of graphs with the largest possible order  for
a given maximum degree and diameter, known as  the $(\Delta,D)$-problem, is a question of particular interest in graph theory.
This problem receives much attention due to its implications  in the design of topologies in interconnection networks 
and its relevance to issues such as  data alignment and the description of  cryptographic protocols. Additionally, the 
 $(\Delta,D)$-problem relates to various graph properties, including regularity, connectivity, and cycle structure.

 Hoffman and Singleton introduced the concept of Moore graphs \cite{HoSi60},   %named after E. F. Moore, 
and proved that  the largest possible order of a graph with maximum degree $\Delta$  ($\Delta > 2$) and diameter $D$  is bounded by 
 \vspace{-3mm}
$$1+\Delta+\Delta(\Delta-1)+\dots+\Delta(\Delta-1)^{D-1}=
    \frac{\Delta(\Delta-1)^{D}-2}{\Delta-2}=n(\Delta,D)\vspace{-3mm}$$
This value is called the {\it Moore bound}, and a graph attaining it is known as a
{\it    Moore graph}. 
Different authors have proved that for $D\geq 2$ and $\Delta\geq 3$, there can only exist Moore graphs for  $D=2$ and $\Delta=3, 7 \mbox{ and } 57$. For the first two cases the graphs are unique and are, respectively, the Petersen's graph on 10 vertices  and the Hoffman--Singleton's graph on 50 vertices. The existence of a Moore graph for  $D=2$  and $\Delta=57$  is not known.
Thus, it is of interest to find graphs which   for a given diameter and maximum degree, have a number of vertices  as close as possible to the Moore bound.

% different approaches
The $(\Delta, D)$ problem for undirected graphs  has been approached in different ways.
It is possible to give bounds to the order of the graphs 
for a given maximum degree and diameter.
On the other hand, as the theoretical bounds are difficult 
to attain, most of the work deals with the construction of  
graphs that,  for this given diameter and maximum degree, have 
a number of vertices as close as possible to the theoretical  bounds.

%different constructive techniques
Different techniques have been developed depending on
the way graphs are generated and their parameters
are calculated.

%Cayley graphs
Roughly half of the  graphs in Table 1 correspond to Cayley
graphs and voltage graphs based on a Cayley graph and have been found  by computer methods \cite{LoSi06,MiSi13,Com25,Com26}. 
%computational methods
The computer is
used for generating the graphs and testing for the desired 
properties. 

%Compound graphs
Compounding is another technique that has proven useful in  producing $(\Delta, D)$  graphs.
This technique  generalizes a method  introduced by Quisquater \cite{Q87}, which involves replacing a single vertex  in a bipartite Moore graph with a complete graph. 
G\'omez, Fiol and Serra  \cite{GoFiSe93}  further  modified this technique to replace several vertices 
of a given graph with another graph or copies of a graph, 
rearranging the edges appropriately (see also \cite{CoGo95,GoFi85,Go09}).
Compound graphs   have been a basic tool in constructing many large 
$(\Delta,D)$ graphs, particularly for small diameters.

Other large graphs have been found as graph products or from particular methods. 
Table 1 gives the current state of the art for degrees up to sixteen and diameter up to ten.

 \section{Recent results  }

Details for the following updates and their adjacency lists, and information and adjacency lists for most graphs in Table 1 with order less than 20000, can be downloaded from \cite{Com26}. 

\subsection{Addition of vertices}

The $(13,3) = 856$ new graph was obtained in 2021 by Vlad Pelekhaty \cite{Pe21} from the graph $Q_8^\prime d$ with degree 13, diameter 3 and 851 vertices described in \cite{Go09} by adding five new vertices and reconnecting several vertices. The resulting graph is regular, has girth 3 and average distance 2.818817.

\subsection{Construction based on a symmetric graph.}

The $(3,8) = 360$ new graph was obtained in 2018 by Jianxiang Chen \cite{Ch18}. 

This graph is derived from the symmetric Foster graph on 144 vertices with diameter 7 and girth 8 by a complete pairing of its edges  as follows:\\   
Let $F$ be the Foster graph F144A and $\sim$ a complete pairing relation on its edges, see \cite{Com26}.
The new  graph, $G$, is constructed as follows:
The vertex set of $G$ is $V(F) \cup E(F)$.
If $v \in V(F)$ and $u \in V(F)$, then they are not connected in $G$.
If $v \in V(F)$ and $u \in E(F)$, then they are connected in $G$ iff $v$ is incident to $u$ in $F$.
If $v \in E(F)$ and $u \in E(F)$, then they are connected in $G$ iff $v \sim u$ by the pairing relation.
The graph $G$ is not a Cayley graph. It is regular with girth 13 and has average distance 6.122563.
 The SageMath program that % finds the pairing 
produces the graph and the adjacency list can be downloaded from \cite{Com26}. \\

\subsection{Semidirect product of cyclic groups}

M. Abas obtained a Cayley graph $(16,2) = 200$ by considering the  semidirect product 
$ (\Z_{10}\times \Z_{10})\rtimes  \Z_{2} $ and generating set 
$X = A\cup B\cup B^{-1}\cup C\cup C^{-1}$, where 
$A = \{(0,0,1)\}$, $ B = \{ (1,0,1),(1,3,1),(1,7,1),(5,0,1),(5,2,1) \}$ and
$ C = \{ (5, 0, 0), (4, 1, 0), (3, 2, 0) \}$, see \cite{Ab17}.

 A Cayley graph $(9,8) = 1697688$ from the  semidirect product  $\Z_{72}\rtimes_{1413}  \Z_{23579}$  with  generating set  $\{ [8,5958],[27,6086], [37,22093], [33,22621], [36,2717]\}$ was produced by A. Rodriguez in 2013, see \cite{Ro13} .

The following table provides details of  the additional new Cayley graph entries in Table 1 obtained by the author from semidirect products of cyclic groups \cite{Com25,Com26}.  
\bigskip

  {\footnotesize \vskip 0.2 cm\hskip -1cm
\renewcommand{\arraystretch}{1.3} 
$\begin{array}{|c|c|c|c|}
\hline\hline
 ( \Delta, D)   &  \mbox{ \bf order}  &  \mbox{ \bf group }      & \mbox{ \bf generators } \\\hline\hline
     (6,8)   &   76\,891    &  \Z_{17}\rtimes_{891}  \Z_{4523}  & \mbox{\footnotesize  [6,1326],[4,1336],[14,1686] } \\ \cdashline{1-4}[2pt/2pt]
   (7,6)   &   12\,264    &  \Z_{24}\rtimes_{90}  \Z_{511}  & \mbox{\footnotesize  [13,77],[6,157],[15,50],[12,7] } \\ \cdashline{1-4}[2pt/2pt]
    (7,7)   &   53\,020    &  \Z_{20}\rtimes_{729}  \Z_{2651}  & \mbox{\footnotesize  [6,894],[17,2271],[18,2411],[10,1210] } \\ \cdashline{1-4}[2pt/2pt]
    (8,5)  &  5\,115       &  \Z_{113}\rtimes_{390}  \Z_{196}  & \mbox{\footnotesize [13,277],[1,290],[4.21],[10,258]. } \\ \cdashline{1-4}[2pt/2pt]       
   (9,4)  & 1\,640  &  \Z_{40}\rtimes_{24}  \Z_{41}  & \mbox{\footnotesize [25,28],[14,40],[29,11],[39,12],[20,35]} \\ \cdashline{1-4}[2pt/2pt]       
 (10,4)  & 2\,331  &  \Z_{9}\rtimes_{44}  \Z_{259}  & \mbox{\footnotesize [8,132],[2,171],[2,71],[4,236],[6,240]}  \\ \cdashline{1-4}[2pt/2pt]
 (10,5) & 13\,203 &  \Z_{81}\rtimes_{22}  \Z_{163}  & \mbox{\footnotesize   [49,70], [64,134], [[78,95], [45,156], [14,90]}  \\ \cdashline{1-4}[2pt/2pt]
 (11,5) & 19\,620 &  \Z_{36}\rtimes_{434}  \Z_{545}  &  \mbox{\footnotesize [22,21], [30,484], [22,513], [33,116 ], [28,421], [18,285] }  \\ \cdashline{1-4}[2pt/2pt]
 (12,5) & 29\,621 &  \Z_{19}\rtimes_{1205}  \Z_{1559}  & \mbox{\footnotesize  [4,358], [15,963], [12,47], [9,233], [14,645], [12,1195]. } \\ \cdashline{1-4}[2pt/2pt]
 (13,5) & 40\,488 &  \Z_{24}\rtimes_{362}  \Z_{1687}  & \mbox{\footnotesize  [1,1454], [5,1427], [2,1659], [15,837], [13,1606], [19,1105], [12,1029].} \\ \cdashline{1-4}[2pt/2pt]
 (14,5) & 58\,095 &  \Z_{45}\rtimes_{191}  \Z_{1291}  & \mbox{\footnotesize  [31,28], [32,290], [28,326], [41,665], [18,278], [24,148], [36,259]. } \\ \cdashline{1-4}[2pt/2pt]
 (15,5) & 77\,520 &  \Z_{48}\rtimes_{772}  \Z_{1615}  & \mbox{\footnotesize  [3,482],[28,1131],[31,682],[47,1424],[2,831],[10,300],[23,1068],[24,0]. }  \\\hline\hline
  \end{array}$
 }

\subsection{A Cayley graph from a group of order $648$}

In a computational search for optimal graphs in the degree--diameter problem for Cayley graphs of order less than $1000$, Marston Conder found a $(5,5)=648$ Cayley graph arising from the $279$th group of order $648$ in the {SmallGroups} database \cite{Con26}. Conder considered a generating set  $S=\{a,b,c,x,y\}$, where $a,b,c$ are involutions  and $x,y$ are elements of order $9$ satisfying $x*y=1$, see \cite{Com26}.

%Finitely presented group G on 5 generators, |G| = 648$ with relations 
%$a:=G.1; b:=G.1 * G.2 * G.7; c:=G.1 * G.5^2 * G.6 * G.7; x:=G.2 * G.3^2; y:=G.2^2 * G.3 * G.4^2 * G.5;$

\section{Table of the largest  $(\Delta,D)$ graphs  }

\bigskip
%%%%%%%%%%%%%%%%%%%%%%%%%%%%%%%%%%%%%%%%%%%%%%%%%%%%%%%%%%%%   Table  Charles Delorme  style                 %%%%%%%%%
%%%%%%%%%%%%%%%%%%%%%%%%%%%%%%%%%%%%%%%%%%%%%%%%%%%%%%%%
\newbox\labox\newbox\nubox\newdimen\lawi\newdimen\nuwi\newdimen\supwi
\def\degre#1{\vbox to 6mm{\vfil\hbox to 6mm{{\large\sl #1}\hfil}}}

\def\NE#1#2{%
  \setbox\labox\hbox{{\it #1}}%
  \setbox\nubox\hbox{{\normalsize
{\bf #2}}}%
  \lawi=\wd\labox\nuwi=\wd\nubox
  \ifdim\nuwi>\lawi\supwi=\nuwi\else\supwi=\lawi\fi
  \vbox to 10mm{% Adjust the height as needed
    \vfil
    \hbox to\supwi{\hfil\box\labox}\vspace{2mm}% Adjust the vertical space between #1 and #2 here
    \hbox to\supwi{\hfil\box\nubox}\vfil
  }
}

\def\BB#1#2{%
  \setbox\labox\hbox{{\it #1}}%
  \setbox\nubox\hbox{#2}%
  \lawi=\wd\labox\nuwi=\wd\nubox
  \ifdim\nuwi>\lawi\supwi=\nuwi\else\supwi=\lawi\fi
  \vbox to 10mm{% Adjust the height as needed
    \vfil
    \hbox to\supwi{\hfil\box\labox}\vspace{2mm}% Adjust the vertical space between #1 and #2 here
    \hbox to\supwi{\hfil\box\nubox}\vfil
  }
}

\def\Bb#1#2{\setbox\labox\hbox{{\it #1}}
\setbox\nubox\hbox{ #2}
\lawi=\wd\labox\nuwi=\wd\nubox
\ifdim\nuwi>\lawi\supwi=\nuwi\else\supwi=\lawi\fi
\vbox to 6mm{\vfil
\hbox to\supwi{\hfil\box\labox}\vfil
\hbox to\supwi{\hfil\box\nubox}\vfil}
}
\def\st{\mathord{*}}

\footnotesize

\renewcommand{\tabcolsep}{0.7mm}
\renewcommand{\arraystretch}{0.9} %1.4 abans
\setlength{\doublerulesep}{2\arrayrulewidth}

\noindent\begin{tabular}{||r||r|r|r|r|r|r|r|r|r||}
\hline\hline
   {\large\sl  D}
 & {\large\sl  2}
 & {\large\sl  3}
 & {\large\sl  4}
 & {\large\sl  5}
 & {\large\sl  6}
 & {\large\sl  7}
 & {\large\sl  8}
& {\large\sl  9}
& {\large\sl  10} \\
 \makebox[6mm][l]{\large$\Delta$} &   &   &   &   &   &  & &  & \\
\hline\hline
\degre{3}&
\BB{P}{10}&\BB{$C_5\st F_4$}{20}&\BB{vC}{38}&\BB{vC}{70}&
\BB{$Exoo$}{132}&\BB{$Exoo$}{196}&\NE{$Chen$}{360}&\BB{$Exoo$}{600}&\BB{$Conder$}{1\,250}\\
\hline
\degre{4}&
\BB{$K_3\st C_5$}{15}&\BB{Allwr}{41}&\BB{$Exoo$}{98}&
\BB{$H'_3$}{364}&\BB{$H_3(K_3)$}{740}&\BB{$Loz$}{1\,320}&\BB{$Loz$}{3\,243}&\BB{$Loz$}{7\,575}&\BB{$Loz$}{17\,703}\\
\hline
\degre{5}&
\BB{$K_3\st X_8$}{24}&\BB{$Exoo$}{72}&\BB{$Exoo$}{212}&
\NE{$Conder$}{648}&\BB{$H_4(K_3)$}{2\,772}&\BB{$Loz$}{5\,516}&
\BB{$Loz$}{17\,030}&\BB{$Loz$}{57\,840}&\BB{$Loz$}{187\,056}\\
\hline
\degre{6}&
\BB{$K_4\st X_8$}{32}&\BB{$Exoo$}{111}&\BB{$Loz$}{390}&
\BB{$Loz$}{1\,404}&\BB{$H_5(K_4)$}{7\,917}&\BB{$Loz$}{19\,383}&
\NE{$Com$}{76\,891}&\BB{$Rod$}{331\,387}&\BB{$Loz$}{1\,253\,615}\\
\hline
\degre{7}&
\BB{HS}{50}&\BB{$Exoo$}{168}&\BB{$Sa$}{672}&
\BB{DH}{2\,756}&\NE{$Com$}{12\,264}&
\NE{$Com$}{53\,020}&\BB{$Loz$}{249\,660}&\BB{$Loz$}{1\,223\,050}&\BB{$Loz$}{6\,007\,230}\\
\hline
\degre{8}&
\BB{$P'_7$}{57}&\BB{$CM,Sa$}{253}&\BB{$Loz$}{1\,100}&
\NE{$Com$}{5\,115}&\BB{$H_7(K_5)$}{39\,672}&\BB{$Loz$}{131\,137}&
\BB{$Loz$}{734\,820}&\BB{$Loz$}{4\,243\,100}&\BB{$Loz$}{24\,897\,161}\\
\hline
\degre{9}&
\BB{$P'_8d$}{74}&\BB{$Q'_8$}{585}&\NE{$Com$}{1\,640}&
\BB{$Rod$}{8\,268}&\BB{$H_8(K_6)$}{75\,893}&\BB{$Loz$}{279\,616}&
\NE{$Rod$}{1\,697\,688}&
\BB{$Loz$}{12\,123\,288}&
\BB{$Loz$}{65\,866\,350}\\
\hline
\degre{10}&
\BB{$P'_9$}{91}&\BB{$Q'_8d$}{650}&\NE{$Com$}{2\,331}&
\NE{$Com$}{13\,203}&\BB{$H_9(K_6)$}{134\,690}&\BB{$Loz$}{583\,083}&
\BB{$Loz$}{4\,293\,452}&
\BB{$Loz$}{27\,997\,191}&
\BB{$Loz$}{201\,038\,922}\\
\hline
\degre{11}&
\BB{$Exoo$}{104}&\BB{$Q'_8d$}{715}&\BB{$Q_7(T_4)$}{3\,200}&
\NE{$Com$}{19\,620}&\BB{$H_7(T_4)$}{156\,864}&
\BB{$Loz$}{1\,001\,268}&
\BB{$Loz$}{7\,442\,328}&
\BB{$Loz$}{72\,933\,102}&
\BB{$Loz$}{600\,380\,000}\\
\hline
\degre{12}&
\BB{$P'_{11}$}{133}&\BB{$Q'_8d^+$}{786}&\BB{$Q'_8\st X_8$}{4\,680}&
\NE{$Com$}{29\,621}&\BB{$H_{11}(K_8)$}{359\,772}&
\BB{$Loz$}{1\,999\,500}&
\BB{$Loz$}{15\,924\,326}&
\BB{$Loz$}{158\,158\,875}&
\BB{$Loz$}{1\,506\,252\,500}\\
\hline
\degre{13}&
\BB{\it MMS}{162}&\NE{$Pel$}{856}&\BB{$Q_9(T_4)$}{6\,560}&
\NE{$Com$}{40\,488}&\BB{$H_9(T_4)$}{531\,440}&
\BB{$Loz$}{3\,322\,080}&
\BB{$Loz$}{29\,927\,790}&
\BB{$Loz$}{249\,155\,760}&
\BB{$Loz$}{3\,077\,200\,700}\\
\hline
\degre{14}&
\BB{$P'_{13}$}{183}&\BB{$Q'_8d^+$}{916}&\BB{$Q_9(T_5)$}{8\,200}&
\NE{$Com$}{58\,095}&\BB{$H_{13}(K_{10})$}{816\,294}&
\BB{$K_1\Sigma_8 H_{11}$}{6\,200\,460}&
\BB{$Loz$}{55\,913\,932}&
\BB{$Loz$}{600\,123\,780}&
\BB{$Loz$}{7\,041\,746\,081}\\
\hline
\degre{15}&  
\BB{$P'_{13}d$}{187}&\BB{$(\otimes Q_{2,4})'$}{1215}&\BB{$Q_{11}(T_4)$}{11\,712}&
\NE{$Com$}{77\,520}& \BB{$H_{11}(T_4)$}{1\,417\,248}&
\BB{$Loz$}{8\,599\,986}&
\BB{$Loz$}{90\,001\,236}&
\BB{$Loz$}{1\,171\,998\,164}&
\BB{$Loz$}{10\,012\,349\,898}\\
\hline
\degre{16}&
\NE{$Abas$}{200}&\BB{$(\otimes Q_3)'$}{1\,600}&\BB{$Q_{11}(T_5)$}{14\,640}&
\BB{$(\otimes H_3)'$}{132\,496}&\BB{$H_{11}(T_5)$}{1\,771\,560}&
\BB{$K_1\Sigma_8 H_{13}$}{14\,882\,658}&
\BB{$Loz$}{140\,559\,416}&\BB{$Loz$}{2\,025\,125\,476}&
\BB{$Loz$}{12\,951\,451\,931}\\
\hline\hline
\end{tabular}
\normalsize
\begin{center}
Table 1: Largest $(\Delta,D)$-graphs (January 2026). In bold new results since the publication of the 2013 survey  \cite{MiSi13}.  Detailed information and adjacency lists for most of the graphs with order less than 20000 can be downloaded from \cite{Com26}.
\end{center}
%%%%%%%%%%%%%%%%%%%%%%%%%%%%%%%%%%%%%%%%%%%%%%%%%%%%
\newpage

\begin{center}
\newdimen\reste
\reste=\textwidth
\advance\reste by -22mm
\long\def\ar#1#2{\noindent\makebox[20mm][l]{#1}\parbox[t]{\reste}{#2}\par}
\bibliographystyle{plain}
{\bf Graphs}\par
\ar{\it Abas}{Graph found by M. Abas \cite{Ab17}}
\ar{\it Allwr}{Special graphs found by J. Allwright \cite{A90}}
\ar{\it Chen}{Graph found by J, Chen \cite{Ch18}}
\ar{\it CM}{Cayley graph, $ \Z_{M}\!\!\rtimes_{A}\!  \Z_{N} $,  found by F. Comellas and M. Mitjana }
\ar{\it Com}{Cayley graphs, $ \Z_{M}\!\!\rtimes_{A}\!  \Z_{N} $, found by F. Comellas (this note) }
\ar{\it Conder}{Graphs found by M. Conder \cite{CoMa13,Con26}}
\ar{{\it vC}}{Compound graphs designed by C. von Conta \cite{vC83}}
\ar{$DH$}{Cayley graph, $ \Z_{M}\!\!\rtimes_{A}\!  \Z_{N} $,  found by M. J. Dinneen and P. Hafner \cite{DH94}}
%\ar{{\it GFS}}{Special graphs $K\Sigma_8 H$  by J. G\'omez, M.A. Fiol and O. Serra \cite{GoFiSe93}}
\ar{$Exoo$}{Graphs built by G. Exoo  \cite{Ex01,Ex10}.}
\ar{$H_q$}{Incidence graph of a regular generalized hexagon \cite{B66}}
\ar{{\it HS}}{Hoffman-Singleton graph}
\ar{$K_n$}{Complete graph}
\ar{$Loz$}{Graphs built by E. Loz  \cite{LoSi06}.}
\ar{\it MMS}{Graph built by B. D. McKay, M. Miller, J. \v{S}ir\'{a}\v{n} \cite{McMiSi98}}
\ar{$P$}{Petersen graph}
\ar{\it Pel}{Graph found by V. Pelekhaty \cite{Pe21}}
\ar{$P_q$}{Incidence graph of projective plane \cite{Hae80}}
\ar{$Q_q$}{Incidence graph of a regular generalized quadrangle \cite{B66}}
\ar{$Rod$}{Cayley graphs,$ \Z_{M}\!\!\rtimes_{A}\!  \Z_{N} $, found by A. Rodriguez \cite{Ro13}}
\ar{$Sa$}{Cayley graphs,$ \Z_{M}\!\!\rtimes_{A}\!  \Z_{N} $, found by M. Sampels \cite{Sa97}}
\ar{$T$}{Tournament}

{\bf Operations}\par
\ar{$G*H$}{Twisted product of graphs \cite{BDF84}}
\ar{$Gd$}{Duplication of some vertices of $G$ \cite{D85a,DF84}}
\ar{$Gd^+$}{Special duplication of vertices of $G$ \cite{Go09}}
\ar{$B'$}{Quotient of the bipartite graph $B$ by a polarity \cite{D85a}}
\ar{$\otimes B$}{The component with polarity of the cartesian product of a bipartite graph $B$ by itself \cite{D85b}}
\ar{$B(K)$}{Substitution of vertices of a bipartite graph $B$ by complete 
     graphs $K$  \cite{CoGo95,Go05,Go09,PGMP09}}
\ar{$B(T)$}{Compound using a bipartite graph $B$ and a tournament $T$
     \cite{GoFi85}}
%\ar{$B(K)<B$}{Compound of $B(K)$ and a bipartite graph $B$ \cite{GFS92}}
\ar{$G\Sigma_iH$}{Various compounding operations \cite{GoFiSe93}}\par
\ar{$ \Z_{M}\!\!\rtimes_{A}\!  \Z_{N} $}{Semidirect product of cyclic groups \cite{DH94}}
%\ar{$DH*$}{Semidirect product of cyclic groups and $Z_n \times Z_n$   \cite{DH92}}
%\ar{$DH**$}{Semidirect product $G \cdot G$ where G is a semidirect product of cyclic groups \cite{DH92}}

\bigskip

Table 2: Graphs and operations of Table 1.
\end{center}

\subsection*{Acknowledgements}

We would like to acknowledge the invaluable contributions of Charles Delorme, who sadly passed away in 2015. The format of this table and detailed information on the graphs and graphs operations are based on his work.

\end{document}